\documentclass[11pt]{article}
\usepackage[utf8]{inputenc}
\usepackage{amsmath,amssymb,amsthm}
\usepackage[a4paper,margin=2.4cm]{geometry}
\usepackage{enumitem}
\setlist[enumerate]{itemsep=5pt}
\newcommand{\dsum}{\displaystyle\sum}
\newcommand{\dprod}{\displaystyle\prod}
\newtheorem{theorem}{Theorem}[section]
\newtheorem{lemma}[theorem]{Lemma}

\newtheorem{corollary}[theorem]{Corollary}
\theoremstyle{remark}
\newtheorem{remark}[theorem]{Remark}

\title{Inverse trigonometric and hyperbolic identities,\\ series and products for Lucas sequences}
\author{Hac\`ene Belbachir\\[4pt]
\small USTHB, Faculty of Mathematics, RECITS Laboratory,\\[-2pt]
\small Po. Box 32 El Alia, 16111 Bab Ezzouar, Algiers, Algeria\\[-2pt]
\small \texttt{hacenebelbachir@gmail.com}, \texttt{hbelbachir@usthb.dz}\\[-2pt]
\small ORCID: 0000-0001-8540-3033}
\date{}

\begin{document}
\maketitle

\begin{abstract}
In 1936 Lehmer observed that the inverse tangents of the reciprocals of the odd-indexed Fibonacci numbers add up exactly to a quarter of $\pi$. Many formulas of the same kind have been found since, each with its own proof: for the inverse tangent and for the inverse hyperbolic tangent, for sums and for differences, for the Fibonacci numbers and for the Lucas numbers. We show that these are not different results but one and the same computation, seen through three signs. A single two-term identity, valid for the generalized Fibonacci sequence and its companion Lucas sequence of any parameter, contains them all. From it follow, with short telescoping proofs, infinitely many series for $\pi$ whose speed of convergence can be chosen at will, series whose sums are multiples of the logarithm of the golden ratio and of its analogues, and infinite and finite products with closed values. The Fibonacci and Lucas cases are listed in full; the classical formulas of Lehmer, of Hoggatt and Ruggles, of Melham and Shannon and of Melham appear among them, and the Pell cases, among them those of Mahon and Horadam, the Jacobsthal cases, where the logarithm of the golden ratio is replaced by that of $2$, the Mersenne numbers and the Chebyshev polynomials follow by substitution.
\end{abstract}



\section{Introduction}

Arctangent series involving Fibonacci numbers go back to Lehmer's Monthly problem \cite{leh}, the evaluation $\sum_{m\geq 1}\tan ^{-1}(1/F_{2m+1})=\pi /4$, and were pursued by Hoggatt and Ruggles \cite{hog}, by Mahon and Horadam \cite{mah} for the Pell numbers, and by Melham and Shannon \cite{melsha} and Melham \cite{mel} for generalized Fibonacci numbers; the latter two works also consider the hyperbolic analogues and the infinite products they produce by exponentiation. The subject was later systematized by Adegoke \cite{ade,ade2} and Frontczak \cite{fro}. In all these works the circular and the hyperbolic identities are proved separately, with the parities of the indices exchanged from one side to the other, the sums and the differences of two inverse tangents are treated as two different identities, and each identity is established for the generalized Fibonacci sequence and then again, separately, for its companion Lucas sequence.

The observation underlying this paper is that these eight computations (two functions, two operations, two sequences) are one and the same. Writing $T_{1}=\tan ^{-1}$ and $T_{-1}=\tanh ^{-1}$, the addition formulas of the two functions read, for signs $\theta =\pm 1$ and $\eta =\pm 1$,
\begin{equation}
T_{\theta }(a)-\eta \,T_{\theta }(b)=T_{\theta }\left( \frac{a-\eta b}{1+\eta \theta ab}\right) . \label{addition}
\end{equation}
Let $W_{1,n}$ be the generalized Fibonacci sequence of parameter $x$ multiplied by $\sqrt{x^{2}+4}$, and $W_{2,n}$ its companion Lucas sequence, both written in the Binet form (\ref{binet}) of Section 2; this normalization removes the factor $\sqrt{x^{2}+4}$ from the identities. Take $a=W_{1,k}/W_{i,n-k}$ and $b=W_{1,k}/W_{i,n+k}$ in (\ref{addition}). The denominator of the right-hand side becomes
\begin{equation*}
1+\eta \theta ab=\frac{W_{i,n+k}W_{i,n-k}+\eta \theta W_{1,k}^{2}}{W_{i,n+k}W_{i,n-k}},
\end{equation*}
and the Catalan identity (\ref{cat}), $W_{i,n+k}W_{i,n-k}=W_{i,n}^{2}+(-1)^{i}(-1)^{n-k}W_{1,k}^{2}$, reduces its numerator to $W_{i,n}^{2}$ as soon as $\eta \theta =(-1)^{i-1}(-1)^{n-k}$, that is for $\theta =(-1)^{n-k+i-1}\eta $; (\ref{addition}) is then Theorem \ref{T1} below. With the three signs $(-1)^{n-k}$, $\eta $ and $(-1)^{i}$, each block of sixteen relations of the author's thesis \cite{these} becomes a single relation, and the case distinctions on the parities of $n$ and $k$, on the nature of the operation and on the sequence disappear into the parity of $n-k$ and the values of $\eta $ and $i$. The paper presents the results of \cite{these} in this form: two-term identities (Section 3), series and their values (Section 4), infinite and finite products (Section 5), a second family built on the duplication formula of the two functions, whose Pell cases appeared in \cite{belben} (Section 6), and the complete list of the Fibonacci and Lucas specializations, with the classical formulas identified among them (Section 7); the Pell and Pell--Lucas numbers, the Jacobsthal and Jacobsthal--Lucas numbers, the Mersenne numbers and the Chebyshev polynomials are reached by substitution (Remark \ref{R-other}), and so are the balancing numbers (Remark \ref{R-bal}). The dictionary for reading an identity involving $T_{\theta }$ is the following: $\theta =1$ gives arctangents and $\theta =-1$ hyperbolic arctangents; $\eta =1$ gives a difference and a telescoping series, $\eta =-1$ a sum and an alternating series.

\section{Preliminaries}

Throughout, $x$ is a positive real number, $\Delta =x^{2}+4$, $\delta =\sqrt{\Delta }$, $\alpha =(x+\delta )/2$, the dominant root of $z^{2}=xz+1$, whose other root is $(x-\delta )/2=-1/\alpha $; the index $i$ takes the values $1$ and $2$, and we put, in the form of the Binet formulas,
\begin{equation}
W_{i,n}=\alpha ^{n}+(-1)^{i+n}\alpha ^{-n}\qquad (n\in \mathbb{Z}). \label{binet}
\end{equation}
In terms of the generalized Fibonacci sequence $U_{n}$ and its companion Lucas sequence $V_{n}$, $U_{0}=0$, $U_{1}=1$, $V_{0}=2$, $V_{1}=x$, $u_{n+1}=xu_{n}+u_{n-1}$ for $u=U$ and $u=V$, one has $W_{1,n}=\delta U_{n}$ and $W_{2,n}=V_{n}$; both $W_{i,n}$ satisfy the same recurrence, $W_{1,0}=0$, $W_{1,1}=\delta $, $W_{2,0}=2$, $W_{2,1}=x$, and $W_{i,-n}=(-1)^{i+n}W_{i,n}$. For $x=1$ one gets $W_{1,n}=\sqrt{5}F_{n}$ and $W_{2,n}=L_{n}$ with the Fibonacci and Lucas numbers, for $x=2$ the Pell and Pell--Lucas numbers $W_{1,n}=2\sqrt{2}P_{n}$, $W_{2,n}=Q_{n}$. Nothing below depends on $x$ being an integer: the identities of Lemma \ref{ids} are identities between the numbers $\alpha ^{n}$, $\alpha ^{-n}$, and the only role of the condition $x>0$ is to ensure $\alpha >1$, hence $W_{i,n}\rightarrow \infty $, which is what the series need. Negative values of $x$ reduce to positive ones through $U_{n}(-x)=(-1)^{n+1}U_{n}(x)$, $V_{n}(-x)=(-1)^{n}V_{n}(x)$. The identities we use are classical (see \cite{hor} or \cite[Chapter 5]{kos}); in the notation (\ref{binet}) they take the following form.

\begin{lemma}
\label{ids}For all integers $n$ and $k$,
\begin{align}
W_{i,n+k}-\eta W_{i,n-k} &=\begin{cases}W_{3-i,n}W_{1,k} & \text{if }\eta =(-1)^{k},\\ W_{i,n}W_{2,k} & \text{if }\eta =-(-1)^{k},\end{cases} \label{add}\\
W_{i,n+k}W_{i,n-k} &=W_{i,n}^{2}+(-1)^{i}(-1)^{n-k}W_{1,k}^{2}, \label{cat}\\
W_{2,n}^{2}-W_{1,n}^{2}=4(-1)^{n},\qquad W_{1,n}W_{2,n} &=W_{1,2n},\qquad W_{2,n}\pm W_{1,n}=2(\pm 1)^{n}\alpha ^{\pm n}. \label{simson}
\end{align}
\end{lemma}

The relations (\ref{simson}) are immediate from (\ref{binet}), and (\ref{add}), (\ref{cat}) are obtained by expanding the products $W_{i,n}W_{j,k}$ and $W_{i,n+k}W_{i,n-k}$ by (\ref{binet}). For each parity of $k$, (\ref{add}) is the usual pair of addition and subtraction formulas.

Besides the addition formula (\ref{addition}), we use the duplication formula
\begin{equation}
2T_{\theta }(t)=T_{\theta }\left( \frac{2t}{1-\theta t^{2}}\right) . \label{dup}
\end{equation}
Both are valid for $ab>-1$ (difference, $\eta =1$), $ab<1$ (sum, $\eta =-1$) and $t^{2}<1$ when $\theta =1$, and for $|a|,|b|,|t|<1$ when $\theta =-1$; they are the two real forms of the single relation $\tanh ^{-1}z=-i\tan ^{-1}(iz)$. In each statement below the sign $\theta $ is specified, and in the proofs $T$ stands for $T_{\theta }$; the identities are stated under the standing assumption that the arguments of $T$ which occur satisfy the above conditions. Since $W_{i,n}\rightarrow \infty $, this holds as soon as the indices are large enough; for $x\geq 1$ the sequences are non-decreasing from $n=1$ on ($U_{1}=U_{2}$ when $x=1$), every $n>k$ is admissible for the arguments of $\tan ^{-1}$ in Section 3, and $n>2k$ for those of $\tanh ^{-1}$, except for $x=1$ and $k=1$, where $n\geq 4$ is needed. The excluded small cases only matter at the start of the series, where they are handled by a uniform rule. Writing $x=2\sinh \tau $ with $\tau >0$ gives $\alpha =e^{\tau }$ and $W_{i,n}=2\cosh n\tau $ or $2\sinh n\tau $ according as $i+n$ is even or odd, so that $\log \alpha =\tau $ and the identities below are, in the end, identities of hyperbolic trigonometry in the variable $\tau $; their interest lies in the arithmetic reading, in which the $W_{i,n}$ are integers or integer multiples of $\delta $.

The sequences $U_{n}(x,y)$, $V_{n}(x,y)$ with recurrence $u_{n+1}=xu_{n}+yu_{n-1}$ and $y>0$ reduce to the case $y=1$: an immediate induction gives $y^{-(n-1)/2}U_{n}(x,y)=U_{n}(x/\sqrt{y},1)$ and $y^{-n/2}V_{n}(x,y)=V_{n}(x/\sqrt{y},1)$, so that every identity below holds for them after the substitutions $U_{n}\rightarrow y^{-(n-1)/2}U_{n}$, $V_{n}\rightarrow y^{-n/2}V_{n}$, $\Delta \rightarrow (x^{2}+4y)/y$. For instance, the Jacobsthal and Jacobsthal--Lucas numbers, which are $U_{n}(1,2)$ and $V_{n}(1,2)$, correspond to the parameter $1/\sqrt{2}$ (Remark \ref{R-other}). The same scaling reduces $y<0$ to $y=-1$, and this case is covered too: for $x>2$ the roots of $z^{2}=xz-1$ are $\alpha ^{\prime 2}$ and $\alpha ^{\prime -2}$, where $\alpha ^{\prime }$ is the dominant root for the parameter $x^{\prime }=\sqrt{x-2}$ and $y=1$, so that $V_{n}(x,-1)=V_{2n}(x^{\prime },1)$ and $x^{\prime }U_{n}(x,-1)=U_{2n}(x^{\prime },1)$: the identities for $y=-1$ are those of the present paper read at even indices. This covers the bisections $F_{2n}=U_{n}(3,-1)$, $L_{2n}=V_{n}(3,-1)$, the Mersenne numbers $2^{n}-1=U_{n}(3,-2)$ and, through $\mathcal{T}_{n}(x)=\frac{1}{2}V_{n}(2x,-1)$ and $\mathcal{U}_{n-1}(x)=U_{n}(2x,-1)$ (see \cite[Remark 7]{belben}), the Chebyshev polynomials $\mathcal{T}_{n}$, $\mathcal{U}_{n}$ of the first and second kind for $x>1$. We shall not repeat these remarks.

\section{Two-term identities}

The identities of this section are the building blocks of all the series and products below. In the first theorem the two arguments have the same numerator $W_{1,k}$ and the sign $\theta $ is imposed by the Catalan identity; in the second, the arguments are quotients of consecutive terms and both signs $\theta $ are allowed. Each theorem is a single relation in $W_{i,n}$: $i=1$ is the identity for the sequence $U$, $i=2$ the identity for the sequence $V$.

\begin{theorem}
\label{T1}For all integers $n>k\geq 1$, with $\theta =(-1)^{n-k+i-1}\eta $,
\begin{equation*}
T_{\theta }\left( \frac{W_{1,k}}{W_{i,n-k}}\right) -\eta \,T_{\theta }\left( \frac{W_{1,k}}{W_{i,n+k}}\right) =T_{\theta }\left( \frac{W_{1,k}\left( W_{i,n+k}-\eta W_{i,n-k}\right) }{W_{i,n}^{2}}\right) .
\end{equation*}
\end{theorem}

\begin{proof}
Put $\varepsilon =(-1)^{n-k}$, so that $\theta =(-1)^{i-1}\eta \varepsilon $ and $\eta \theta =(-1)^{i-1}\varepsilon $. By (\ref{addition}) with $a=W_{1,k}/W_{i,n-k}$ and $b=W_{1,k}/W_{i,n+k}$,
$$T(a)-\eta \,T(b)=T\left( \frac{W_{1,k}\left( W_{i,n+k}-\eta W_{i,n-k}\right) }{W_{i,n+k}W_{i,n-k}+(-1)^{i-1}\varepsilon W_{1,k}^{2}}\right) ,$$
and the denominator is $W_{i,n}^{2}$ by (\ref{cat}), since $\varepsilon =(-1)^{n-k}$.
\end{proof}

By (\ref{add}), the right-hand side of Theorem \ref{T1} is
\begin{equation}
T_{\theta }\Big( \frac{W_{1,k}^{2}\,W_{3-i,n}}{W_{i,n}^{2}}\Big) \quad \text{if }\eta =(-1)^{k},\qquad T_{\theta }\Big( \frac{W_{1,k}W_{2,k}}{W_{i,n}}\Big) =T_{\theta }\Big( \frac{W_{1,2k}}{W_{i,n}}\Big) \quad \text{if }\eta =-(-1)^{k}. \label{T1explicit}
\end{equation}
For $x=1$ and $k=1,2$ these are the sixteen relations of Corollary \ref{C1}.

\begin{theorem}
\label{T2}For all integers $n>k\geq 1$ and for both signs $\theta =\pm 1$,
\begin{equation*}
T_{\theta }\left( \frac{W_{i,n-k}}{W_{i,n}}\right) -\eta \,T_{\theta }\left( \frac{W_{i,n}}{W_{i,n+k}}\right) =T_{\theta }\left( \frac{(1-\eta )W_{i,n}^{2}+(-1)^{n-k+i}W_{1,k}^{2}}{W_{i,n}\left( W_{i,n+k}+\theta \eta W_{i,n-k}\right) }\right) .
\end{equation*}
\end{theorem}

\begin{proof}
Put $\varepsilon =(-1)^{n-k}$. With $a=W_{i,n-k}/W_{i,n}$ and $b=W_{i,n}/W_{i,n+k}$ one has
\begin{equation*}
a-\eta b=\frac{W_{i,n-k}W_{i,n+k}-\eta W_{i,n}^{2}}{W_{i,n}W_{i,n+k}}=\frac{(1-\eta )W_{i,n}^{2}+(-1)^{i}\varepsilon W_{1,k}^{2}}{W_{i,n}W_{i,n+k}},\qquad 1+\eta \theta ab=\frac{W_{i,n+k}+\theta \eta W_{i,n-k}}{W_{i,n+k}},
\end{equation*}
the first by (\ref{cat}); (\ref{addition}) concludes.
\end{proof}

The denominator is given by (\ref{add}): $W_{i,n+k}+\theta \eta W_{i,n-k}$ equals $W_{i,n}W_{2,k}$ if $\theta \eta =(-1)^{k}$ and $W_{3-i,n}W_{1,k}$ if $\theta \eta =-(-1)^{k}$. The six theorems of \cite{melsha} are the cases $k=1$, $\eta =1$ of Theorems \ref{T1} and \ref{T2}, read with $\tan ^{-1}$ and with $\tanh ^{-1}$.

\section{Series}

The relation of Theorem \ref{T1} telescopes ($\eta =1$) or alternates ($\eta =-1$) along an arithmetic progression $n=2sm+r$ of step $2s$, because then $n+s=(2m+1)s+r$ is the value of $n-s$ for $m+1$; the parity of $n-s$, which is that of $s+r$, does not depend on $m$, so that the sign $\theta $ of Theorem \ref{T1} is the same for all the terms, and the factor $\eta ^{m-1}$ below is $1$ or $(-1)^{m-1}$ according to the case.

\begin{theorem}
\label{T3}Let $s\geq 1$, $0\leq r<2s$ and $\theta =(-1)^{s+r+i-1}\eta $. Then
\begin{equation*}
\sum_{m\geq 1}\eta ^{m-1}\,T_{\theta }\left( \frac{W_{1,s}\left( W_{i,(2m+1)s+r}-\eta W_{i,(2m-1)s+r}\right) }{W_{i,2sm+r}^{2}}\right) =T_{\theta }\left( \frac{W_{1,s}}{W_{i,s+r}}\right) .
\end{equation*}
When the first terms are not defined, the series starts at the first index $m_{0}$ for which all terms are defined, and the right-hand side becomes $\eta ^{m_{0}-1}$ times $T_{\theta }\big( W_{1,s}/W_{i,(2m_{0}-1)s+r}\big) $. For $x\geq 1$ this happens only for $r=0$, with $i=1$ and $\eta =-(-1)^{s}$ or with $i=2$, $s$ odd and $\eta =-1$ (the argument of $\tanh ^{-1}$ on the right being then $\geq 1$), and, for $x=1$ only, for $(s,r)=(1,1)$ with $i=1$ and $\eta =-1$, where the first term is $\tanh ^{-1}\left( L_{3}/F_{3}^{2}\right) =\tanh ^{-1}1$; in all these cases $m_{0}=2$.
\end{theorem}

\begin{proof}
Put $a_{l}=T\left( W_{1,s}/W_{i,l}\right) $ and apply Theorem \ref{T1} with $k=s$, $n=2sm+r$: the $m$-th term is $a_{(2m-1)s+r}-\eta \,a_{(2m+1)s+r}$, so that
\begin{equation*}
\sum_{m=1}^{N}\eta ^{m-1}\left( a_{(2m-1)s+r}-\eta \,a_{(2m+1)s+r}\right) =a_{s+r}-\eta ^{N}a_{2sN+s+r}\rightarrow a_{s+r}.
\end{equation*}
When the series starts at $m_{0}$, the surviving term is $\eta ^{m_{0}-1}a_{(2m_{0}-1)s+r}$.
\end{proof}

The right-hand side takes two different forms according to the function. On the circular side, $\tan ^{-1}\left( W_{1,s}/W_{i,s+r}\right) =\frac{\pi }{2}-\tan ^{-1}\left( W_{i,s+r}/W_{1,s}\right) $, and for $i=1$, $r=0$ this is $\pi /4$. On the hyperbolic side and for $i=2$, $r=0$, $s$ even and $\eta =1$, (\ref{simson}) gives $\tanh ^{-1}\left( W_{1,s}/W_{2,s}\right) =\frac{1}{2}\log \frac{W_{2,s}+W_{1,s}}{W_{2,s}-W_{1,s}}=s\log \alpha $. We single out these cases.

\begin{corollary}
\label{C-pi}For every $s\geq 1$,
\begin{gather}
\sum_{m\geq 1}\tan ^{-1}\frac{U_{2s}^{2}\,V_{4sm}}{U_{4sm}^{2}}=\sum_{m\geq 1}\left( -1\right) ^{m-1}\tan ^{-1}\frac{U_{2s-1}^{2}\,V_{(4s-2)m}}{U_{(4s-2)m}^{2}}=\frac{\pi }{4}, \label{pi4}\\
\sum_{m\geq 1}\tanh ^{-1}\frac{\delta ^{3}U_{2s}^{2}\,U_{4sm}}{V_{4sm}^{2}}=2s\log \alpha . \label{logalpha}
\end{gather}
\end{corollary}

\begin{proof}
Take $r=0$ in Theorem \ref{T3} with $(i,s,\eta )\rightarrow (1,2s,1)$, $(1,2s-1,-1)$ and $(2,2s,1)$, use (\ref{add}) for the numerators and $W_{1,n}=\delta U_{n}$, $W_{2,n}=V_{n}$.
\end{proof}

The two series for $\pi /4$ and the series for $\log \alpha $ are thus the three readings of one telescoping; what is $\pi /4$ on the circular side is a multiple of $\log \alpha $ on the hyperbolic side. In the parametrization $x=2\sinh \tau $ of Section 2 the third series equals $2s\tau $ and the corresponding product of Theorem \ref{T4} equals $e^{4s\tau }$; for $\tau =1$, that is $x=e-e^{-1}$, one has $\alpha =e$, $\delta =2\cosh 1$, $U_{n}=\sinh n/\cosh 1$ or $\cosh n/\cosh 1$ and $V_{n}=2\cosh n$ or $2\sinh n$ according as $n$ is even or odd, and for instance
\begin{equation*}
\sum_{m\geq 1}\tanh ^{-1}\frac{2\sinh ^{2}2\,\sinh 4m}{\cosh ^{2}4m}=2,\qquad \prod_{m\geq 1}\frac{V_{4m}^{2}+\delta ^{3}U_{2}^{2}U_{4m}}{V_{4m}^{2}-\delta ^{3}U_{2}^{2}U_{4m}}=e^{4}.
\end{equation*}
The number $e$ enters here only through the parameter: for $x=1$ the base is $\varphi $, and the telescoping method produces values which are algebraic numbers, logarithms of units or rational multiples of $\pi $, never $e$ itself. Since $U_{n}\sim \alpha ^{n}/\delta $, the $m$-th term of the first series is asymptotic to $\Delta U_{2s}^{2}\alpha ^{-4sm}$, so that the convergence, geometric with ratio $\alpha ^{-4s}$, is as fast as one wishes, and the remainder after $N$ terms is explicit, namely $\tan ^{-1}\left( U_{2s}/U_{4sN+2s}\right) $.

\section{Products}

The relation $e^{2\tanh ^{-1}a}=\frac{1+a}{1-a}$ turns every identity read with $\tanh ^{-1}$ into a product. The infinite products come from Theorem \ref{T3}, the finite ones from Theorem \ref{T2}.

\begin{theorem}
\label{T4}Let $s\geq 1$ and $0\leq r<2s$. If $(-1)^{s+r+i-1}\eta =-1$,
\begin{equation*}
\prod_{m\geq 1}\left( \frac{W_{i,2sm+r}^{2}+W_{1,s}\left( W_{i,(2m+1)s+r}-\eta W_{i,(2m-1)s+r}\right) }{W_{i,2sm+r}^{2}-W_{1,s}\left( W_{i,(2m+1)s+r}-\eta W_{i,(2m-1)s+r}\right) }\right) ^{\eta ^{m-1}}=\frac{W_{i,s+r}+W_{1,s}}{W_{i,s+r}-W_{1,s}},
\end{equation*}
with the starting rule of Theorem \ref{T3}: for $x\geq 1$ and $r=0$, the product with $i=1$ and $\eta =-(-1)^{s}$ starts at $m=2$ and has the value $\left( \frac{W_{1,3s}+W_{1,s}}{W_{1,3s}-W_{1,s}}\right) ^{\eta }$, the product with $i=2$, $s$ odd and $\eta =-1$ starts at $m=2$ and has the value $\frac{W_{2,3s}-W_{1,s}}{W_{2,3s}+W_{1,s}}$; for $x=1$ the product with $(i,s,r,\eta )=(1,1,1,-1)$ also starts at $m=2$, with the value $\frac{F_{4}-F_{1}}{F_{4}+F_{1}}=\frac{1}{2}$.
\end{theorem}

\begin{proof}
Exponentiate Theorem \ref{T3} in the cases where $\theta =(-1)^{s+r+i-1}\eta =-1$.
\end{proof}

For $i=2$, $r=0$, $s$ even and $\eta =1$, the value of the product is $\frac{W_{2,s}+W_{1,s}}{W_{2,s}-W_{1,s}}=\alpha ^{2s}$ by (\ref{simson}): these infinite products are powers of $\alpha $, hence units of $\mathbb{Q}(\delta )$ when $x$ is an integer. The finite products below have, by contrast, values which are quotients of terms of the sequences.

\begin{theorem}
\label{T5}Let $s\geq 1$, $r\geq 0$, $N\geq 1$ and put
\begin{equation*}
Z_{i,m}=\frac{(1-\eta )W_{i,sm+r}^{2}+(-1)^{s(m-1)+r+i}W_{1,s}^{2}}{W_{i,sm+r}\left( W_{i,s(m+1)+r}-\eta W_{i,s(m-1)+r}\right) }.
\end{equation*}
Then, provided $W_{i,ls+r}<W_{i,(l+1)s+r}$ for $l=0,\dots ,N$,
\begin{equation*}
\prod_{m=1}^{N}\left( \frac{1+Z_{i,m}}{1-Z_{i,m}}\right) ^{\eta ^{m-1}}=\frac{W_{i,s+r}+W_{i,r}}{W_{i,s+r}-W_{i,r}}\cdot \left( \frac{W_{i,sN+s+r}-\eta W_{i,sN+r}}{W_{i,sN+s+r}+\eta W_{i,sN+r}}\right) ^{\eta ^{N-1}}.
\end{equation*}
For $x\geq 1$ the condition excludes only $(x,s,r)=(1,1,0)$ and $(1,1,1)$ when $i=1$, and $s=1$, $r=0$ with $1\leq x\leq 2$ when $i=2$, where $W_{2,0}=2\geq W_{2,1}=x$.
\end{theorem}

\begin{proof}
Put $b_{n}=\tanh ^{-1}\left( W_{i,n-s}/W_{i,n}\right) $; then $\tanh ^{-1}\left( W_{i,n}/W_{i,n+s}\right) =b_{n+s}$, and $e^{2b_{n}}$ equals $\left( W_{i,n}+W_{i,n-s}\right) /\left( W_{i,n}-W_{i,n-s}\right) $. Theorem \ref{T2} with $\theta =-1$, $k=s$ and $n=sm+r$ reads $b_{sm+r}-\eta \,b_{s(m+1)+r}$$=\tanh ^{-1}Z_{i,m}$; multiplying by $\eta ^{m-1}$ and summing over $m=1,\dots ,N$ gives $b_{s+r}-\eta ^{N}b_{sN+s+r}$, and one exponentiates. The condition on the indices is that $b_{ls+r}$ be defined; for $x\geq 1$ it fails only in the cases listed, where $\tanh ^{-1}(F_{1}/F_{2})$, resp.\ $\tanh ^{-1}(V_{0}/V_{1})=\tanh ^{-1}(2/x)$, is not defined.
\end{proof}

By (\ref{add}), for $\eta =1$ the factor $Z_{i,m}$ equals
\begin{equation*}
(-1)^{s(m-1)+r+i}\frac{W_{1,s}}{W_{1,2sm+2r}}\ (s\text{ even}),\qquad (-1)^{s(m-1)+r+i}\frac{W_{1,s}^{2}}{W_{2,s}W_{i,sm+r}^{2}}\ (s\text{ odd});
\end{equation*} for even $s$ the two relations $i=1,2$ with $\eta =1$ are therefore equivalent to each other, since (\ref{add}) gives, for every $n$,
\begin{equation*}
\frac{W_{1,s+n}+W_{1,n}}{W_{1,s+n}-W_{1,n}}\cdot \frac{W_{2,s+n}+W_{2,n}}{W_{2,s+n}-W_{2,n}}=\frac{W_{2,s}+2}{W_{2,s}-2}.
\end{equation*} Since $W_{i,sN+r}/W_{i,sN+s+r}\rightarrow \alpha ^{-s}$, letting $N\rightarrow \infty $ for $\eta =1$ gives the infinite products
\begin{equation*}
\prod_{m\geq 1}\frac{1+Z_{i,m}}{1-Z_{i,m}}=\frac{W_{i,s+r}+W_{i,r}}{W_{i,s+r}-W_{i,r}}\cdot \frac{\alpha ^{s}-1}{\alpha ^{s}+1},
\end{equation*}
in which $\frac{\alpha ^{s}-1}{\alpha ^{s}+1}=\tanh \frac{s\log \alpha }{2}$; for $\eta =-1$ the partial products have no limit: those of odd and of even rank $N$ tend respectively to $A\,\frac{\alpha ^{s}+1}{\alpha ^{s}-1}$ and $A\,\frac{\alpha ^{s}-1}{\alpha ^{s}+1}$, where $A=\frac{W_{i,s+r}+W_{i,r}}{W_{i,s+r}-W_{i,r}}$.

\section{Identities built on the duplication formula}

A different choice of numerators leads to a second family. Instead of two arguments with the same numerator $W_{1,k}$, one takes two arguments with the same \emph{constant} numerator $2$; this is the only constant for which $W_{i,n+k}W_{i,n-k}\pm 4$ is a sum or a difference of two squares of terms, $W_{i,n}^{2}\pm W_{j,k}^{2}$ or $W_{3-i,n}^{2}\pm W_{j,k}^{2}$, by (\ref{cat}) and (\ref{simson}), and the duplication formula (\ref{dup}) then produces a factor $2$. For the Pell numbers, where $V_{1}=2$, the arguments $2/W_{2,n}=Q_{1}/Q_{n}$ are quotients of two Pell--Lucas numbers, which is how this family first arose in \cite{belben}: the Pell case of the relation below, with $i=2$, $k=1$, $\eta =1$ and $n$ even, is the two-term identity of relation 3 of Theorem 3 of \cite{belben}. The sequence which appears on the right depends on the sign $\eta (-1)^{k}$: for $k\geq 1$ we let
\begin{equation}
j=1\ \text{if }\eta (-1)^{k}=1,\qquad j=2\ \text{if }\eta (-1)^{k}=-1, \label{WZ}
\end{equation}
and $j^{\prime }=j$ if $i=2$, $j^{\prime }=3-j$ if $i=1$: the quotient $W_{j,k}/W_{j^{\prime },n}$ below is $U_{k}/U_{n}$ or $V_{k}/V_{n}$ for $i=2$, and $\delta U_{k}/V_{n}$ or $V_{k}/\left( \delta U_{n}\right) $ for $i=1$.

\begin{theorem}
\label{T6}For all integers $n>k\geq 1$, with $\theta =(-1)^{n-k+i-1}$ and $j,j^{\prime }$ as in (\ref{WZ}),
\begin{equation*}
T_{\theta }\left( \frac{2}{W_{i,n-k}}\right) -\eta \,T_{\theta }\left( \frac{2}{W_{i,n+k}}\right) =2\,T_{\theta }\left( \frac{W_{j,k}}{W_{j^{\prime },n}}\right) .
\end{equation*}
\end{theorem}

\begin{proof}
Put $\varepsilon =(-1)^{n-k}$, so that $\theta =(-1)^{i-1}\varepsilon $. By (\ref{addition}) the left-hand side is
\begin{equation*}
T\left( \frac{2\left( W_{i,n+k}-\eta W_{i,n-k}\right) }{W_{i,n+k}W_{i,n-k}+4(-1)^{i-1}\eta \varepsilon }\right) ,
\end{equation*}
and by (\ref{cat}) the denominator is $W_{i,n}^{2}+(-1)^{i}\varepsilon \left( W_{1,k}^{2}-4\eta \right) $. If $j=2$, that is $\eta =-(-1)^{k}$, then $W_{1,k}^{2}-4\eta =W_{1,k}^{2}+4(-1)^{k}=W_{2,k}^{2}$ by (\ref{simson}), the denominator is $W_{i,n}^{2}+(-1)^{i}\varepsilon W_{2,k}^{2}$ and the numerator is $2W_{i,n}W_{2,k}$ by (\ref{add}): the argument is $\frac{2t}{1-(-1)^{i-1}\varepsilon t^{2}}$ with $t=W_{2,k}/W_{i,n}$, and $j^{\prime }=i$ in this case. If $j=1$, that is $\eta =(-1)^{k}$, then $(-1)^{n}=(-1)^{n-k}(-1)^{k}=\varepsilon \eta $, so that $W_{i,n}^{2}-4(-1)^{i}\varepsilon \eta =W_{i,n}^{2}-4(-1)^{i}(-1)^{n}=W_{3-i,n}^{2}$ by (\ref{simson}): the denominator is $W_{3-i,n}^{2}+(-1)^{i}\varepsilon W_{1,k}^{2}$, the numerator is $2W_{3-i,n}W_{1,k}$, the argument is $\frac{2t}{1-(-1)^{i-1}\varepsilon t^{2}}$ with $t=W_{1,k}/W_{3-i,n}$, and $j^{\prime }=3-i$. The duplication formula (\ref{dup}) concludes.
\end{proof}

Along $n=2sm+r$ one has $n+s=(2m+1)s+r$, the value of $n-s$ for $m+1$, and the relation telescopes for $\eta =1$ and alternates for $\eta =-1$.

\begin{theorem}
\label{T7}Let $s\geq 1$, $0\leq r<2s$, $\theta =(-1)^{s+r+i-1}$ and $j,j^{\prime }$ as in (\ref{WZ}) for $k=s$. Then
\begin{equation*}
\sum_{m\geq 1}\eta ^{m-1}\,T_{\theta }\left( \frac{W_{j,s}}{W_{j^{\prime },2sm+r}}\right) =\frac{1}{2}T_{\theta }\left( \frac{2}{W_{i,s+r}}\right) .
\end{equation*}
\end{theorem}

\begin{proof}
With $g(l)=T\left( 2/W_{i,l}\right) $, Theorem \ref{T6} shows that the $m$-th term is $\frac{1}{2}\big( g( (2m-1)s+r) -\eta \,g( (2m+1)s+r) \big) $; the partial sums are therefore $\frac{1}{2}\big( g(s+r)-\eta ^{N}g(2sN+s+r)\big) $, which tend to $\frac{1}{2}g(s+r)$. For $x\geq 1$ all terms are defined, since $W_{1,l}>2$ for $l\geq 1$, $W_{2,l}>2$ for $l\geq 2$, for $i=2$, $s+r=1$ one has $\theta =1$, and the arguments $W_{j,s}/W_{j^{\prime },2sm+r}$ are at most $1$, with equality only for $x=1$, $s=1$, $i=2$, $\eta =-1$, $m=1$, where $\theta =1$; for $0<x<1$ the starting rule of Theorem \ref{T3} applies.
\end{proof}

On the circular side the value is $\frac{\pi }{4}-\frac{1}{2}\tan ^{-1}\frac{W_{i,s+r}}{2}$; it equals $\pi /8$ exactly when $W_{i,s+r}=2$, that is, for $x\geq 1$, for $i=2$, $s=1$, $r=0$ and $x=2$ (the Pell numbers), and, for $0<x<1$, for $i=1$, $s+r=2$ and $x^{2}=2\sqrt{2}-2$. On the hyperbolic side, exponentiation gives products whose values are the square roots
\begin{equation*}
\rho _{i,l}=\sqrt{\frac{W_{i,l}+2}{W_{i,l}-2}},
\end{equation*}
and since $W_{2,l}^{2}=W_{2,2l}+2(-1)^{l}$ and $W_{1,l}^{2}=W_{2,2l}-2(-1)^{l}$ by (\ref{binet}), $\rho _{2,2l}=W_{2,l}/W_{1,l}$ for $l$ even and $W_{1,l}/W_{2,l}$ for $l$ odd.

\begin{corollary}
\label{C-dup}Let $s\geq 1$, $0\leq r<2s$, $N\geq 1$ and $j,j^{\prime }$ as in (\ref{WZ}) for $k=s$. If $(-1)^{s+r+i-1}=-1$,
\begin{equation*}
\prod_{m=1}^{N}\left( \frac{W_{j^{\prime },2sm+r}+W_{j,s}}{W_{j^{\prime },2sm+r}-W_{j,s}}\right) ^{\eta ^{m-1}}=\rho _{i,s+r}\,\rho _{i,2sN+s+r}^{-\eta ^{N}},
\end{equation*}
and as $N\rightarrow \infty $ the product tends to $\rho _{i,s+r}$.
\end{corollary}

\begin{proof}
Exponentiate the partial sums in the proof of Theorem \ref{T7}, using $e^{g(l)}=\rho _{i,l}$ when the function is $\tanh ^{-1}$, that is $(-1)^{s+r+i-1}=-1$.
\end{proof}

The constant $2$ also relates the quotients of Theorem \ref{T2} for the two sequences; relation 1 of Theorem 3 of \cite{belben} is the Pell case of the first relation below.

\begin{theorem}
\label{T8}For all integers $n>k\geq 1$, with $\theta =(-1)^{k+1}$,
\begin{enumerate}
\item $T_{\theta }\left( \left( \dfrac{W_{1,n\pm k}}{W_{1,n}}\right) ^{\mp 1}\right) +T_{\theta }\left( \left( \dfrac{W_{2,n\pm k}}{W_{2,n}}\right) ^{\mp 1}\right) =T_{\theta }\left( \dfrac{2}{W_{2,k}}\right) ,$

\item $T_{\theta }\left( \left( \dfrac{W_{1,n\pm k}}{W_{1,n}}\right) ^{\mp 1}\right) -T_{\theta }\left( \left( \dfrac{W_{2,n\pm k}}{W_{2,n}}\right) ^{\mp 1}\right) =\left( -1\right) ^{n+1}\,T_{\theta }\left( \dfrac{2}{W_{2,k\pm 2n}}\right) .$
\end{enumerate}
Consequently,
\begin{equation}
T_{\theta }\left( \frac{W_{i,n-k}}{W_{i,n}}\right) =\frac{1}{2}\left( T_{\theta }\frac{2}{W_{2,k}}+(-1)^{n+k+i}\,T_{\theta }\frac{2}{W_{2,2n-k}}\right) . \label{dup1}
\end{equation}
\end{theorem}

\begin{proof}
Put $\varepsilon =(-1)^{k}$, so that $\theta =-\varepsilon $. By (\ref{addition}), the argument of $T$ in the sum of relation 1 is
\begin{equation*}
\frac{W_{1,n-k}W_{2,n}+W_{2,n-k}W_{1,n}}{W_{1,n}W_{2,n}+\varepsilon W_{1,n-k}W_{2,n-k}}=\frac{2W_{1,2n-k}}{W_{1,2n}+\varepsilon W_{1,2n-2k}}=\frac{2}{W_{2,k}},
\end{equation*}
the first equality by (\ref{simson}) and $W_{1,m}W_{2,n}+W_{2,m}W_{1,n}=2W_{1,m+n}$, the second by (\ref{add}) with $\eta =-(-1)^{k}$, which gives $W_{1,2n}+(-1)^{k}W_{1,2n-2k}=W_{1,2n-k}W_{2,k}$. In the difference the argument is
\begin{equation*}
\frac{W_{1,n-k}W_{2,n}-W_{2,n-k}W_{1,n}}{W_{1,n}W_{2,n}-\varepsilon W_{1,n-k}W_{2,n-k}}=\frac{2(-1)^{n+k+1}W_{1,k}}{W_{1,2n}-(-1)^{k}W_{1,2n-2k}}=\frac{2(-1)^{n+k+1}}{W_{2,2n-k}},
\end{equation*}
since $W_{1,m}W_{2,n}-W_{2,m}W_{1,n}=2(-1)^{n}W_{1,m-n}$ and $W_{1,-k}=(-1)^{k+1}W_{1,k}$, and by (\ref{add}) with $\eta =(-1)^{k}$; here $(-1)^{n+k+1}=-\varepsilon (-1)^{n}$. The versions with $n+k$ follow by replacing $n$ by $n+k$, which changes $(-1)^{n}$ into $(-1)^{n}(-1)^{k}=\varepsilon (-1)^{n}$. Since $W_{2,-l}=(-1)^{l}W_{2,l}$ by (\ref{binet}), $(-1)^{n+k+1}/W_{2,2n-k}=(-1)^{n+1}/W_{2,k-2n}$. Half-sum and half-difference give (\ref{dup1}).\end{proof}

\section{The Fibonacci and Lucas numbers}

The specialization $x=1$ is collected here in full, together with the places where the classical formulas appear; the other sequences are obtained by substitution, as explained in Remark \ref{R-other}. For $x=1$ one has $\delta =\sqrt{5}$, $\alpha =\varphi =(1+\sqrt{5})/2$, $W_{1,n}=\sqrt{5}F_{n}$, $W_{2,n}=L_{n}$, $F_{1}=F_{2}=L_{1}=1$ and $L_{2}=3$; in every quotient of two terms $W_{1,m}/W_{1,n}$ the factor $\sqrt{5}$ cancels. Each relation below whose sign $\theta $ depends on the parity of $n$ contains one circular and one hyperbolic identity, and each relation with a sign $\eta $ contains a sum and a difference.

\subsection{Two-term identities}

The sixteen relations of Theorem \ref{T1} with $k=1,2$ and $i=1,2$, written with (\ref{T1explicit}), are the following.

\begin{corollary}
\label{C1}In the $\tan ^{-1}$ readings for all $n\geq 2$ in relations 1 and 3 and all $n\geq 3$ in relations 2 and 4, in the $\tanh ^{-1}$ readings for all $n\geq 5$ and for $n=3,4$ whenever the arguments of $\tanh ^{-1}$ are less than $1$, with $\theta =(-1)^{n+1}\eta $ in relations 1 and 4 and $\theta =(-1)^{n}\eta $ in relations 2 and 3,
\begin{enumerate}
\item $T_{\theta }\left( \dfrac{1}{F_{n-1}}\right) -\eta \,T_{\theta }\left( \dfrac{1}{F_{n+1}}\right) =T_{\theta }\left( \dfrac{1}{F_{n}}\right) \ (\eta =1),\quad T_{\theta }\left( \dfrac{L_{n}}{F_{n}^{2}}\right) \ (\eta =-1),$

\item $T_{\theta }\left( \dfrac{1}{F_{n-2}}\right) -\eta \,T_{\theta }\left( \dfrac{1}{F_{n+2}}\right) =T_{\theta }\left( \dfrac{L_{n}}{F_{n}^{2}}\right) \ (\eta =1),\quad T_{\theta }\left( \dfrac{3}{F_{n}}\right) \ (\eta =-1),$

\item $T_{\theta }\left( \dfrac{\sqrt{5}}{L_{n-1}}\right) -\eta \,T_{\theta }\left( \dfrac{\sqrt{5}}{L_{n+1}}\right) =T_{\theta }\left( \dfrac{\sqrt{5}}{L_{n}}\right) \ (\eta =1),\quad T_{\theta }\left( \dfrac{5\sqrt{5}F_{n}}{L_{n}^{2}}\right) \ (\eta =-1),$

\item $T_{\theta }\left( \dfrac{\sqrt{5}}{L_{n-2}}\right) -\eta \,T_{\theta }\left( \dfrac{\sqrt{5}}{L_{n+2}}\right) =T_{\theta }\left( \dfrac{5\sqrt{5}F_{n}}{L_{n}^{2}}\right) \ (\eta =1),\quad T_{\theta }\left( \dfrac{3\sqrt{5}}{L_{n}}\right) \ (\eta =-1).$
\end{enumerate}
\end{corollary}

Several of these relations, in their circular or hyperbolic readings, appear in \cite{hog} and \cite{mel}. The sixteen relations of Theorem \ref{T2} with $k=1,2$ and $i=1,2$ are the following; the denominators $F_{n+k}+\theta \eta F_{n-k}$ and $L_{n+k}+\theta \eta L_{n-k}$ of Theorem \ref{T2}, which do not depend on the parity of $n$, are given by (\ref{add}) as displayed.

\begin{corollary}
\label{C2}For all $n\geq 2$ in relations 1 and 2 and $n\geq 3$ in relations 3 and 4, except $n=2$ in relation 1 with $\theta =-1$, and for both signs $\theta =\pm 1$,
\begin{enumerate}
\item $T_{\theta }\left( \dfrac{F_{n-1}}{F_{n}}\right) -\eta \,T_{\theta }\left( \dfrac{F_{n}}{F_{n+1}}\right) =T_{\theta }\left( \dfrac{(1-\eta )F_{n}^{2}+(-1)^{n}}{F_{n}L_{n}}\right) \ (\theta \eta =1),$ \\ $\qquad T_{\theta }\left( \dfrac{(1-\eta )F_{n}^{2}+(-1)^{n}}{F_{n}^{2}}\right) \ (\theta \eta =-1),$

\item $T_{\theta }\left( \dfrac{L_{n-1}}{L_{n}}\right) -\eta \,T_{\theta }\left( \dfrac{L_{n}}{L_{n+1}}\right) =T_{\theta }\left( \dfrac{(1-\eta )L_{n}^{2}-5(-1)^{n}}{5F_{n}L_{n}}\right) \ (\theta \eta =1),$ \\ $\qquad T_{\theta }\left( \dfrac{(1-\eta )L_{n}^{2}-5(-1)^{n}}{L_{n}^{2}}\right) \ (\theta \eta =-1),$

\item $T_{\theta }\left( \dfrac{F_{n-2}}{F_{n}}\right) -\eta \,T_{\theta }\left( \dfrac{F_{n}}{F_{n+2}}\right) =T_{\theta }\left( \dfrac{(1-\eta )F_{n}^{2}-(-1)^{n}}{3F_{n}^{2}}\right) \ (\theta \eta =1),$ \\ $\qquad T_{\theta }\left( \dfrac{(1-\eta )F_{n}^{2}-(-1)^{n}}{F_{n}L_{n}}\right) \ (\theta \eta =-1),$

\item $T_{\theta }\left( \dfrac{L_{n-2}}{L_{n}}\right) -\eta \,T_{\theta }\left( \dfrac{L_{n}}{L_{n+2}}\right) =T_{\theta }\left( \dfrac{(1-\eta )L_{n}^{2}+5(-1)^{n}}{3L_{n}^{2}}\right) \ (\theta \eta =1),$ \\ $\qquad T_{\theta }\left( \dfrac{(1-\eta )L_{n}^{2}+5(-1)^{n}}{5F_{n}L_{n}}\right) \ (\theta \eta =-1).$
\end{enumerate}
\end{corollary}

For instance, relation 1 reads
\begin{gather*}
\tan ^{-1}\tfrac{F_{n-1}}{F_{n}}-\tan ^{-1}\tfrac{F_{n}}{F_{n+1}}=\tan ^{-1}\tfrac{(-1)^{n}}{F_{n}L_{n}},\quad \tanh ^{-1}\tfrac{F_{n-1}}{F_{n}}-\tanh ^{-1}\tfrac{F_{n}}{F_{n+1}}=\tanh ^{-1}\tfrac{(-1)^{n}}{F_{n}^{2}},\\
\tan ^{-1}\tfrac{F_{n-1}}{F_{n}}+\tan ^{-1}\tfrac{F_{n}}{F_{n+1}}=\tan ^{-1}\tfrac{2F_{n}^{2}+(-1)^{n}}{F_{n}^{2}},\quad \tanh ^{-1}\tfrac{F_{n-1}}{F_{n}}+\tanh ^{-1}\tfrac{F_{n}}{F_{n+1}}=\tanh ^{-1}\tfrac{2F_{n}^{2}+(-1)^{n}}{F_{n}L_{n}}.
\end{gather*}

\subsection{Series}

The twenty-four series of Theorem \ref{T3} with $s=1,2$ specialize as follows; the exceptional starting index $m=2$ occurs in the cases foreseen in Theorem \ref{T3}.

\begin{corollary}
\label{C3}For the Fibonacci numbers ($i=1$),
\begin{enumerate}
\item $\dsum \limits_{m\geq 1}\tan ^{-1}\dfrac{1}{F_{2m+1}}=\dsum \limits_{m\geq 1}\tan ^{-1}\dfrac{L_{4m}}{F_{4m}^{2}}=\dsum \limits_{m\geq 1}\left( -1\right) ^{m-1}\tan ^{-1}\dfrac{L_{2m}}{F_{2m}^{2}}=\dfrac{\pi }{4},$

\item $\dsum \limits_{m\geq 1}\tan ^{-1}\dfrac{L_{4m+2}}{F_{4m+2}^{2}}=\dfrac{\pi }{2}-\tan ^{-1}3,$ \\ $\qquad \dsum \limits_{m\geq 1}\left( -1\right) ^{m-1}\tan ^{-1}\dfrac{3}{F_{4m+1}}=\tan ^{-1}\dfrac{1}{2},$ \\ $\qquad \dsum \limits_{m\geq 1}\left( -1\right) ^{m-1}\tan ^{-1}\dfrac{3}{F_{4m+3}}=\tan ^{-1}\dfrac{1}{5},$

\item $\dsum \limits_{m\geq 2}\tanh ^{-1}\dfrac{1}{F_{2m}}=\dsum \limits_{m\geq 1}\tanh ^{-1}\dfrac{L_{4m+1}}{F_{4m+1}^{2}}=\tanh ^{-1}\dfrac{1}{2},$ \\ $\qquad \dsum \limits_{m\geq 1}\tanh ^{-1}\dfrac{L_{4m+3}}{F_{4m+3}^{2}}=\tanh ^{-1}\dfrac{1}{5},$

\item $\dsum \limits_{m\geq 2}\left( -1\right) ^{m-1}\tanh ^{-1}\dfrac{L_{2m+1}}{F_{2m+1}^{2}}=-\tanh ^{-1}\dfrac{1}{3},$ \\ $\qquad \dsum \limits_{m\geq 1}\left( -1\right) ^{m-1}\tanh ^{-1}\dfrac{3}{F_{4m+2}}=\tanh ^{-1}\dfrac{1}{3},$ \\ $\qquad \dsum \limits_{m\geq 2}\left( -1\right) ^{m-1}\tanh ^{-1}\dfrac{3}{F_{4m}}=-\tanh ^{-1}\dfrac{1}{8},$
\end{enumerate}
and for the Lucas numbers ($i=2$),
\begin{enumerate}[start=5]
\item $\dsum \limits_{m\geq 1}\tan ^{-1}\dfrac{\sqrt{5}}{L_{2m}}=\dfrac{\pi }{2}-\tan ^{-1}\dfrac{1}{\sqrt{5}},$ \\ $\qquad \dsum \limits_{m\geq 1}\tan ^{-1}\dfrac{5\sqrt{5}F_{4m+1}}{L_{4m+1}^{2}}=\dfrac{\pi }{2}-\tan ^{-1}\dfrac{4}{\sqrt{5}},$ \\ $\qquad \dsum \limits_{m\geq 1}\tan ^{-1}\dfrac{5\sqrt{5}F_{4m+3}}{L_{4m+3}^{2}}=\dfrac{\pi }{2}-\tan ^{-1}\dfrac{11}{\sqrt{5}},$

\item $\dsum \limits_{m\geq 1}\left( -1\right) ^{m-1}\tan ^{-1}\dfrac{5\sqrt{5}F_{2m+1}}{L_{2m+1}^{2}}=\dsum \limits_{m\geq 1}\left( -1\right) ^{m-1}\tan ^{-1}\dfrac{3\sqrt{5}}{L_{4m}}=\tan ^{-1}\dfrac{\sqrt{5}}{3},$ \\ $\qquad \dsum \limits_{m\geq 1}\left( -1\right) ^{m-1}\tan ^{-1}\dfrac{3\sqrt{5}}{L_{4m+2}}=\tan ^{-1}\dfrac{\sqrt{5}}{7},$

\item $\dsum \limits_{m\geq 1}\tanh ^{-1}\dfrac{\sqrt{5}}{L_{2m+1}}=\dsum \limits_{m\geq 1}\tanh ^{-1}\dfrac{5\sqrt{5}F_{4m}}{L_{4m}^{2}}=\tanh ^{-1}\dfrac{\sqrt{5}}{3}=2\log \varphi ,$ \\ $\qquad \dsum \limits_{m\geq 1}\tanh ^{-1}\dfrac{5\sqrt{5}F_{4m+2}}{L_{4m+2}^{2}}=\tanh ^{-1}\dfrac{\sqrt{5}}{7},$ \\ (the second of these is (\ref{logalpha}) with $s=1$),

\item $\dsum \limits_{m\geq 2}\left( -1\right) ^{m-1}\tanh ^{-1}\dfrac{5\sqrt{5}F_{2m}}{L_{2m}^{2}}=-\tanh ^{-1}\dfrac{\sqrt{5}}{4},$ \\ $\qquad \dsum \limits_{m\geq 1}\left( -1\right) ^{m-1}\tanh ^{-1}\dfrac{3\sqrt{5}}{L_{4m+1}}=\tanh ^{-1}\dfrac{\sqrt{5}}{4},$ \\ $\qquad \dsum \limits_{m\geq 1}\left( -1\right) ^{m-1}\tanh ^{-1}\dfrac{3\sqrt{5}}{L_{4m+3}}=\tanh ^{-1}\dfrac{\sqrt{5}}{11}.$
\end{enumerate}
\end{corollary}

The first series is Lehmer's \cite{leh} (also in \cite{hog,melsha}); relations (2.12)--(2.15) of \cite{mel} are among the series of Corollary \ref{C3}, and relations (2.16), (2.17) of \cite{mel} among the products of Corollary \ref{C4}. The twenty-four series of Theorem \ref{T7} with $s=1,2$ are the following.

\begin{corollary}
\label{C5}For $i=2$ (constant numerator $2$ over the Lucas numbers),
\begin{enumerate}
\item $\dsum \limits_{m\geq 1}\tan ^{-1}\dfrac{1}{L_{2m}}=\dsum \limits_{m\geq 1}(-1)^{m-1}\tan ^{-1}\dfrac{1}{F_{2m}}=\dfrac{1}{2}\tan ^{-1}2,$ \\ $\qquad \dsum \limits_{m\geq 1}\tanh ^{-1}\dfrac{1}{L_{2m+1}}=\dsum \limits_{m\geq 1}(-1)^{m-1}\tanh ^{-1}\dfrac{1}{F_{2m+1}}=\dfrac{1}{4}\log 5,$

\item $\dsum \limits_{m\geq 1}\tanh ^{-1}\dfrac{1}{F_{4m}}=\dsum \limits_{m\geq 1}(-1)^{m-1}\tanh ^{-1}\dfrac{3}{L_{4m}}=\dfrac{1}{4}\log 5,$ \\ $\qquad \dsum \limits_{m\geq 1}\tan ^{-1}\dfrac{1}{F_{4m+1}}=\dsum \limits_{m\geq 1}(-1)^{m-1}\tan ^{-1}\dfrac{3}{L_{4m+1}}=\dfrac{1}{2}\tan ^{-1}\dfrac{1}{2},$

\item $\dsum \limits_{m\geq 1}\tanh ^{-1}\dfrac{1}{F_{4m+2}}=\dsum \limits_{m\geq 1}(-1)^{m-1}\tanh ^{-1}\dfrac{3}{L_{4m+2}}=\dfrac{1}{2}\tanh ^{-1}\dfrac{2}{7},$ \\ $\qquad \dsum \limits_{m\geq 1}\tan ^{-1}\dfrac{1}{F_{4m+3}}=\dsum \limits_{m\geq 1}(-1)^{m-1}\tan ^{-1}\dfrac{3}{L_{4m+3}}=\dfrac{1}{2}\tan ^{-1}\dfrac{2}{11},$
\end{enumerate}
and for $i=1$ (constant numerator $2$ over $\sqrt{5}F_{n}$),
\begin{enumerate}[start=4]
\item $\dsum \limits_{m\geq 1}\tanh ^{-1}\dfrac{1}{\sqrt{5}F_{2m}}=\dsum \limits_{m\geq 1}(-1)^{m-1}\tanh ^{-1}\dfrac{\sqrt{5}}{L_{2m}}=\dfrac{3}{2}\log \varphi ,$ \\ $\qquad \dsum \limits_{m\geq 1}\tan ^{-1}\dfrac{1}{\sqrt{5}F_{2m+1}}=\dsum \limits_{m\geq 1}(-1)^{m-1}\tan ^{-1}\dfrac{\sqrt{5}}{L_{2m+1}}=\dfrac{1}{2}\tan ^{-1}\dfrac{2}{\sqrt{5}},$

\item $\dsum \limits_{m\geq 1}\tan ^{-1}\dfrac{\sqrt{5}}{L_{4m}}=\dsum \limits_{m\geq 1}(-1)^{m-1}\tan ^{-1}\dfrac{3}{\sqrt{5}F_{4m}}=\dfrac{1}{2}\tan ^{-1}\dfrac{2}{\sqrt{5}},$ \\ $\qquad \dsum \limits_{m\geq 1}\tanh ^{-1}\dfrac{\sqrt{5}}{L_{4m+1}}=\dsum \limits_{m\geq 1}(-1)^{m-1}\tanh ^{-1}\dfrac{3}{\sqrt{5}F_{4m+1}}=\dfrac{1}{2}\log \varphi ,$

\item $\dsum \limits_{m\geq 1}\tan ^{-1}\dfrac{\sqrt{5}}{L_{4m+2}}=\dsum \limits_{m\geq 1}(-1)^{m-1}\tan ^{-1}\dfrac{3}{\sqrt{5}F_{4m+2}}=\dfrac{1}{2}\tan ^{-1}\dfrac{2}{3\sqrt{5}},$ \\ $\qquad \dsum \limits_{m\geq 1}\tanh ^{-1}\dfrac{\sqrt{5}}{L_{4m+3}}=\dsum \limits_{m\geq 1}(-1)^{m-1}\tanh ^{-1}\dfrac{3}{\sqrt{5}F_{4m+3}}=\dfrac{1}{2}\tanh ^{-1}\dfrac{2}{5\sqrt{5}}.$
\end{enumerate}
\end{corollary}

\subsection{Products and the identities of Theorem \ref{T8}}

The products of Theorems \ref{T4} and \ref{T5} and of Corollary \ref{C-dup} are the exponentials of the hyperbolic series above; those with the simplest values, and the identities of Theorem \ref{T8}, are the following.

\begin{corollary}
\label{C4}For the Fibonacci and Lucas numbers,
\begin{enumerate}
\item $\dprod \limits_{m\geq 1}\dfrac{F_{4m+1}^{2}+L_{4m+1}}{F_{4m+1}^{2}-L_{4m+1}}=3,\qquad \dprod \limits_{m\geq 2}\dfrac{F_{6m}+8}{F_{6m}-8}=\dfrac{9}{8},$ \\ $\qquad \dprod \limits_{m\geq 1}\dfrac{F_{6m+1}^{2}+4(-1)^{m-1}L_{6m+1}}{F_{6m+1}^{2}-4(-1)^{m-1}L_{6m+1}}=5,\qquad \dprod \limits_{m\geq 2}\dfrac{F_{4m}+3(-1)^{m-1}}{F_{4m}-3(-1)^{m-1}}=\dfrac{7}{9},$

\item $\dprod \limits_{m\geq 1}\dfrac{L_{4m}^{2}+5\sqrt{5}F_{4m}}{L_{4m}^{2}-5\sqrt{5}F_{4m}}=\varphi ^{4},\qquad \dprod \limits_{m\geq 2}\dfrac{L_{6m}^{2}+20\sqrt{5}(-1)^{m-1}F_{6m}}{L_{6m}^{2}-20\sqrt{5}(-1)^{m-1}F_{6m}}=\dfrac{38-\sqrt{5}}{38+\sqrt{5}},$

\item $\dprod \limits_{m\geq 1}\dfrac{L_{6m+1}+8\sqrt{5}}{L_{6m+1}-8\sqrt{5}}=\dfrac{7+2\sqrt{5}}{7-2\sqrt{5}},\qquad \dprod \limits_{m\geq 1}\dfrac{L_{4m+1}+3\sqrt{5}(-1)^{m-1}}{L_{4m+1}-3\sqrt{5}(-1)^{m-1}}=\dfrac{4+\sqrt{5}}{4-\sqrt{5}},$

\item $\dprod \limits_{m\geq 1}\dfrac{F_{4m}+1}{F_{4m}-1}=\dprod \limits_{m\geq 1}\dfrac{L_{2m+1}+1}{L_{2m+1}-1}=\sqrt{5},\qquad \dprod \limits_{m\geq 1}\dfrac{\sqrt{5}F_{2m}+1}{\sqrt{5}F_{2m}-1}=\varphi ^{3},$ \\ $\qquad \dprod \limits_{m\geq 1}\dfrac{L_{4m+1}+\sqrt{5}}{L_{4m+1}-\sqrt{5}}=\varphi ,\qquad \dprod \limits_{m=1}^{N}\dfrac{F_{4m}+1}{F_{4m}-1}=\dfrac{L_{2N+1}}{F_{2N+1}},$

\item $\tan ^{-1}\dfrac{F_{n-1}}{F_{n}}=\dfrac{1}{2}\left( \tan ^{-1}2+(-1)^{n}\tan ^{-1}\dfrac{2}{L_{2n-1}}\right) ,$ \\ $\qquad \tan ^{-1}\dfrac{L_{n-1}}{L_{n}}=\dfrac{1}{2}\left( \tan ^{-1}2-(-1)^{n}\tan ^{-1}\dfrac{2}{L_{2n-1}}\right) ,$ \\ $\qquad \tan ^{-1}\dfrac{F_{n-1}}{F_{n}}+\tan ^{-1}\dfrac{L_{n-1}}{L_{n}}=\tan ^{-1}2,\qquad \tanh ^{-1}\dfrac{F_{n-2}}{F_{n}}+\tanh ^{-1}\dfrac{L_{n-2}}{L_{n}}=\tanh ^{-1}\dfrac{2}{3}\quad (n\geq 3).$
\end{enumerate}
\end{corollary}

\begin{remark}
\label{R-other}Every relation of this section has an analogue for each value of the parameter, obtained by the substitutions $F_{n}\rightarrow U_{n}$, $L_{n}\rightarrow V_{n}$, $\sqrt{5}\rightarrow \delta $, $5\rightarrow \Delta $, $\varphi \rightarrow \alpha $ and $1,1,1,3\rightarrow U_{1},U_{2},V_{1},V_{2}$, that is, by reading Theorems \ref{T1} to \ref{T8} for the value of $x$ at hand; only the starting indices, governed by the rule of Theorem \ref{T3}, and the numerical values of the sums have to be recomputed. For the Pell and Pell--Lucas numbers ($x=2$, $\delta =2\sqrt{2}$, $\alpha =1+\sqrt{2}$, $U_{1}=1$, $U_{2}=V_{1}=2$, $V_{2}=6$) one obtains in this way, for instance, the dual $\sum_{m\geq 0}\tan ^{-1}\left( 2/P_{2m+1}\right) =\pi /2$ of Lehmer's formula (from Theorem \ref{T3} written with (\ref{T1explicit}), $s=r=1$; it is the series of relation 2 of Theorem 3 of \cite{belben}, and the Pell polynomials are treated in \cite{mah}), the series $\sum_{m\geq 0}\tan ^{-1}\left( 2/Q_{2m}\right) =3\pi /8$ (the series of relation 3 of Theorem 3 of \cite{belben}) and $\sum_{m\geq 1}(-1)^{m-1}\tan ^{-1}\left( 1/P_{2m}\right) =\pi /8$ from Theorem \ref{T7}, both consequences of the Pell case of Theorem \ref{T6} quoted in Section 6, $\tan ^{-1}\frac{P_{n-1}}{P_{n}}+\tan ^{-1}\frac{Q_{n-1}}{Q_{n}}=\frac{\pi }{4}$ from Theorem \ref{T8}, the Pell relations of \cite{mel}, and the products $\prod_{m\geq 1}\frac{P_{4m}+2}{P_{4m}-2}=\prod_{m\geq 1}\frac{Q_{2m+1}+2}{Q_{2m+1}-2}=\sqrt{2}$, $\prod_{m\geq 1}\frac{\sqrt{2}P_{2m}+1}{\sqrt{2}P_{2m}-1}=1+\sqrt{2}$ and $\prod_{m\geq 1}\frac{Q_{4m}^{2}+64\sqrt{2}P_{4m}}{Q_{4m}^{2}-64\sqrt{2}P_{4m}}=(1+\sqrt{2})^{4}$. The sequences with a second parameter $y$ are covered by the reduction of Section 2, and so are the Chebyshev polynomials $\mathcal{T}_{n}(x)$, $\mathcal{U}_{n}(x)$: for $x>1$, reading Corollary \ref{C-pi} for the parameter $\sqrt{2x-2}$, with $V_{2n}=2\mathcal{T}_{n}(x)$ and $U_{2n}=\sqrt{2x-2}\,\mathcal{U}_{n-1}(x)$, gives
\begin{equation*}
\sum_{m\geq 1}\tan ^{-1}\frac{2\mathcal{U}_{s-1}(x)^{2}\,\mathcal{T}_{2sm}(x)}{\mathcal{U}_{2sm-1}(x)^{2}}=\frac{\pi }{4},\quad \sum_{m\geq 1}\tanh ^{-1}\frac{2(x^{2}-1)^{3/2}\,\mathcal{U}_{s-1}(x)^{2}\,\mathcal{U}_{2sm-1}(x)}{\mathcal{T}_{2sm}(x)^{2}}=s\cosh ^{-1}x,
\end{equation*}
since the dominant root is then $x+\sqrt{x^{2}-1}=e^{\cosh ^{-1}x}$. In the same way, the Jacobsthal numbers $J_{n}$ and the Jacobsthal--Lucas numbers $j_{n}$, defined by $J_{0}=0$, $J_{1}=1$, $j_{0}=2$, $j_{1}=1$ and $u_{n+1}=u_{n}+2u_{n-1}$, are $U_{n}(1,2)$ and $V_{n}(1,2)$; the reduction of Section 2 brings them to the parameter $1/\sqrt{2}$, whose dominant root is $\sqrt{2}$, and Corollary \ref{C-pi} gives, for every $s\geq 1$,
\begin{equation*}
\sum_{m\geq 1}\tan ^{-1}\frac{4^{s(m-1)}J_{2s}^{2}\,j_{4sm}}{J_{4sm}^{2}}=\frac{\pi }{4},\qquad \sum_{m\geq 1}\tanh ^{-1}\frac{27\cdot 4^{s(m-1)}J_{2s}^{2}\,J_{4sm}}{j_{4sm}^{2}}=s\log 2,
\end{equation*}
and, for $s=1$, $\sum_{m\geq 1}(-1)^{m-1}\tan ^{-1}\left( 2^{m-1}j_{2m}/J_{2m}^{2}\right) =\pi /4$. With a negative second parameter, the Mersenne numbers $M_{n}=2^{n}-1$ and their companions $2^{n}+1$, which are $U_{n}(3,-2)$ and $V_{n}(3,-2)$, are brought by the same reduction to $y=-1$, and Theorem \ref{T3} gives, for every $s\geq 1$,
\begin{equation*}
\sum_{m\geq 1}\tan ^{-1}\frac{2^{s(m-1)}M_{s}^{2}\,(2^{2sm}+1)}{M_{2sm}^{2}}=\frac{\pi }{4},\qquad \text{for instance}\qquad \sum_{m\geq 1}\tan ^{-1}\frac{2^{m-1}(4^{m}+1)}{(4^{m}-1)^{2}}=\frac{\pi }{4},
\end{equation*}
and, on the hyperbolic side, $\sum_{m\geq 1}(-1)^{m-1}\tanh ^{-1}\left( 15\cdot 2^{2m-1}/M_{4m+2}\right) =\tanh ^{-1}\frac{2}{5}=\frac{1}{2}\log \frac{7}{3}$.
\end{remark}

\begin{remark}
\label{R-bal}The balancing numbers $B_{n}$ and the Lucas-balancing numbers $C_{n}$ of Behera and Panda \cite{behpan}, defined by $B_{0}=0$, $B_{1}=1$, $C_{0}=1$, $C_{1}=3$ and $u_{n+1}=6u_{n}-u_{n-1}$, are $B_{n}=U_{n}(6,-1)=P_{2n}/2$ and $C_{n}=\frac{1}{2}V_{n}(6,-1)=Q_{2n}/2$: they are the Pell and Pell--Lucas numbers read at even indices, as in Section 2. Corollary \ref{C-pi} and Theorem \ref{T3} for the Pell numbers then give, for every $s\geq 1$,
\begin{equation*}
\sum_{m\geq 1}\tan ^{-1}\frac{2B_{s}^{2}\,C_{2sm}}{B_{2sm}^{2}}=\frac{\pi }{4},\qquad \sum_{m\geq 1}(-1)^{m-1}\tanh ^{-1}\frac{B_{2s}}{B_{2sm+r}}=\tanh ^{-1}\frac{B_{s}}{B_{s+r}}\quad (1\leq r<2s),
\end{equation*}
for instance $\sum_{m\geq 1}\tan ^{-1}\left( 2C_{2m}/B_{2m}^{2}\right) =\pi /4$ and $\sum_{m\geq 1}(-1)^{m-1}\tanh ^{-1}\left( 6/B_{2m+1}\right) =\tanh ^{-1}\frac{1}{6}$. Balancing problems in the binomial coefficient triangle are studied in \cite{belsza}.
\end{remark}

\medskip
\noindent\textbf{Declaration.} The mathematical results are the author's. The generative AI assistant Claude (Anthropic) was used to check the identities numerically to high precision and to assist with the typesetting and the English of the manuscript. The author reviewed and takes full responsibility for the content.

\end{document}